\documentclass{amsart}
\usepackage{amssymb}

\usepackage{amscd}

\newtheorem{theorem}{Theorem}[section]

\newtheorem{lemma}[theorem]{Lemma}

\theoremstyle{definition}
\newtheorem{definition}{Definition}[section]

\theoremstyle{remark}

\begin{document}
\title{An analogue of Abel's theorem}
\author{Herbert Clemens}
\address{Mathematics Department, Ohio State University\\
231 W. 18th Ave.\\
Columbus, OH 43210, USA}
\email{clemens@math.ohio-state.edu}
\date{November 18, 2002}
\maketitle

\begin{abstract}
This work makes a parallel construction for curves on threefolds to a
``current-theoretic'' proof of Abel's theorem giving the rational
equivalence of divisors P and Q on a Riemann surface when Q - P is
(equivalent to) zero in the Jacobian variety of the Riemann surface. The
parallel construction is made for homologous ''sub-canonical'' curves P and
Q on a general class of threefolds. If P and Q are algebraically equivalent
and Q - P is zero in the (intermediate) Jacobian of a threefold, the
construction ''almost'' gives rational equivalence.
\end{abstract}

\section{Introduction\label{1}\protect\footnote{%
Partially supported by NSF grant DMS-0200895}}

This work establishes a parallel between

1) a proof of the classical theorem of Abel deriving rational equivalence
classes of divisors $P$ and $Q$ of degree $d$ on a Riemann surface $X$ from
equality of their images in the Jacobian variety $J\left( X\right) $,

and

2) a construction for certain cohomologous curves $P$ and $Q$ on a threefold 
$X$ for which the image of $Q-P$ in the (intermediate) Jacobian $J\left(
X\right) $ is zero.

This paper was motivated by work of Richard Thomas (see \cite{C}). Indeed
the point of view is in large measure due to him. To explain the analogy
with Abel's theorem, one must recast the classical proof of that theorem in
the language of forms with values in distributions, that is, currents (see 
\cite{K}). In that language, the classical proof of Abel's theorem goes
something like this.

Take two effective divisors $P$ and $Q$ of degree $d$ on a Riemann surface $X
$. Consider a one-current 
\begin{equation*}
\Gamma 
\end{equation*}
such that 
\begin{equation*}
Q-P=\partial \left( \Gamma \right) 
\end{equation*}
that is, the operator which assigns to a $C^{\infty }$-one-form on $X$ its
integral over some path $\Gamma $ from $P$ to $Q$. This current is the
pull-back of a current on the multiplicative group $\Bbb{C}^{*}$ of complex
numbers as follows:

For the unique topological line bundle $L_{\infty }$ of degree $d$ on $X$,
there are two complex structures, one giving 
\begin{equation*}
\mathcal{O}\left( P\right)
\end{equation*}
and the other giving 
\begin{equation*}
\mathcal{O}\left( Q\right) .
\end{equation*}
Thus the holomorphic sections $s_{P}$ and $s_{Q}$ with respect to these
complex structures give two different $C^{\infty }$-sections of $L_{\infty }$
which we will continue to call $s_{P}$ and $s_{Q}$. The quotient 
\begin{equation*}
g:=\frac{s_{Q}}{s_{P}}
\end{equation*}
then gives a $C^{\infty }$-map 
\begin{equation*}
X-\left( \left| P\right| \cup \left| Q\right| \right) \rightarrow \Bbb{C}%
^{*}.
\end{equation*}
The current 
\begin{equation*}
\Gamma
\end{equation*}
is just the pull-back of the current on $\Bbb{C}^{*}$ given by the positive
real axis, that is 
\begin{equation*}
\Gamma =g^{*}\left( \left( 0,+\infty \right) \right) .
\end{equation*}

Furthermore, on the compactification $\Bbb{P}^{1}$ of $\Bbb{C}^{*}$ we have
the cohomological relation 
\begin{equation*}
d\mathrm{\log }z\sim 2\pi i\cdot \left( 0,+\infty \right) 
\end{equation*}
of currents. Pulling back via $g^{*}$, we have the equality 
\begin{equation*}
\int\nolimits_{\Gamma }=\frac{1}{2\pi i}\int\nolimits_{X}\alpha \wedge 
\end{equation*}
as functionals on $H^{1,0}\left( X\right) $, where 
\begin{equation*}
\alpha :=d\mathrm{\log }g.
\end{equation*}
Also we note that the $\left( 0,1\right) $-component of $d\mathrm{\log }z$,
that is, the linear operator obtained by integrating $\left( 1,0\right) $%
-forms against $d\mathrm{\log }z$, is actually bounded, that is, is given by
integrating against a $C^{\infty }$-form $\alpha ^{0,1}$ on $X$ of type $%
\left( 0,1\right) $.

Now if $P$ is Jacobian-equivalent to $Q$, the $\left( 0,1\right) $-summand
of the current $\alpha ^{0,1}$ (appropriately normalized by the $\left(
0,1\right) $-summand of an integral cocycle) is 
\begin{equation*}
\overline{\partial }b
\end{equation*}
for some $C^{\infty }$-function $b$. So in this case, 
\begin{equation*}
\psi :=g^{*}\left( d\mathrm{\log }z\right) -db
\end{equation*}
is a form of type $\left( 1,0\right) $ which is $d$-closed on $X-\left(
\left| P\right| \cup \left| Q\right| \right) .$ So $\psi $ is $\overline{%
\partial }$-closed and of type $\left( 1,0\right) $ and therefore
meromorphic, having poles with residues which are integral multiples of $2%
\pi i.$ We complete the proof of Abel's theorem by remarking that 
\begin{equation*}
f=e^{\int \psi }
\end{equation*}
is the rational function giving the rational equivalence of $P$ and $Q$.

After first giving in detail the distribution-theoretic proof of classical
Abel's theorem, we will make an analogous construction in the case in which $%
P$ and $Q$ are certain types of effective algebraic one-cycles on a
threefold $X$. We first produce a rank-2 vector bundle $E_{\infty }$ whose
first Chern class is trivial and whose second Chern class is represented by $%
P$ or by $Q$. We note that a choice of a metric $\mu $ on $E_{\infty }$ is
equivalent to giving $E_{\infty }$ the structure of a quaternionic line
bundle. Then we mimic, for the quaternionic line bundle $E_{\infty }$, all
but the ``trivial'' last step of the proof of Abel's theorem for divisors on
curves. Namely for each metric $\mu $ on $E_{\infty }$ we use Chern-Simons
theory to produce a three-form $\alpha _{\mu ,PQ}$ on $X$ giving the normal
function $Q-P$ and such that:

1) If $P$ and $Q$ are algebraically equivalent, then $\alpha _{\mu
,PQ}^{\left( 1,2\right) +\left( 0,3\right) }$ is $\overline{\partial }$
-closed as a form on $X$.

2) If in addition $P$ and $Q$ are Abel-Jacobi equivalent, then $\alpha _{\mu
,PQ}^{\left( 1,2\right) +\left( 0,3\right) }$ (normalized by the $\left(
1,2\right) +\left( 0,3\right) $ summand of an integral cycle) is $\overline{%
\partial }$-exact as a form on $X$; so there is an associated form $\psi
_{\mu }$ on $X$ of type $\left( 3,0\right) +\left( 2,1\right) $ such that 
\begin{equation*}
\psi _{\mu }-\alpha _{\mu ,PQ}\in dA_{X}^{2}.
\end{equation*}

\noindent Recall that $P$ and $Q$ are \textit{algebraically equivalent} if
they are homologous on some (possibly reducible) divisor on $X$.

Finally, for the quaternionic line bundle $E_{\infty }$ we will use
Chern-Simons theory to produce a $3$-current $\alpha _{PQ}$ giving the
normal function $Q-P$ such that $\alpha _{PQ}$ is $d$-closed on $X^{\prime
}:=X-\left( \left| P\right| \cup \left| Q\right| \right) $ and:

1) If $P$ and $Q$ are algebraically equivalent, then $\alpha _{PQ}^{\left(
1,2\right) +\left( 0,3\right) }$ is $\overline{\partial }$-closed as a
current on $X$.

2) If in addition $P$ and $Q$ are Abel-Jacobi equivalent, then $\alpha
_{PQ}^{\left( 1,2\right) +\left( 0,3\right) }$ (normalized by the $\left(
1,2\right) +\left( 0,3\right) $ summand of an integral cycle) is $\overline{%
\partial }$-exact as a current on $X$; so there is a canonically associated
current $\psi $ on $X$ of type $\left( 3,0\right) +\left( 2,1\right) $ such
that 
\begin{equation*}
\psi -\alpha _{PQ}\in d\left\{ 2-currents\ on\ X\right\} .
\end{equation*}
Furthermore 
\begin{equation*}
\left. \psi \right| _{X^{\prime }}
\end{equation*}
is a $d$-closed $3$-current, that is, its integral against the coboundary of
compactly supported $2$-forms on $X^{\prime }$ is zero.

To understand a potential significance of $\psi $, suppose that $\left. \psi
\right| _{X^{\prime }}$ turns out to be a $3$-form and the algebraic
equivalence of $P$ and $Q$ is given by an algebraic family 
\begin{equation*}
\begin{array}{lll}
S & \overset{s}{\longrightarrow } & X \\ 
\downarrow ^{r} &  &  \\ 
C &  & 
\end{array}
\end{equation*}
such that

1) $C$ is a smooth irreducible curve,

2) $S$ is a smooth surface proper and flat over $C$,

3) for two points $p$ and $q$ in $C$ and the corresponding fibers $S_{p}$
and $S_{q}$ of $r$, we have a rational equivalence 
\begin{equation*}
s_{*}\left( S_{q}\right) -s_{*}\left( S_{p}\right) \equiv Q-P,
\end{equation*}

4) 
\begin{equation*}
s^{-1}\left( \left| P\right| \cup \left| Q\right| \right) \subseteq
S_{p}\cup S_{q}.
\end{equation*}

Then 
\begin{equation*}
r_{*}s^{*}\psi
\end{equation*}
is a $d$-closed form of type $\left( 1,0\right) $ on $C^{\prime }:=C-\left\{
p,q\right\} $ and 
\begin{equation*}
f=e^{\int r_{*}s^{*}\psi }
\end{equation*}
would therefore be a meromorphic function on $C$ giving the rational
equivalence of $P$ and $Q$.

\section{Classical theorem of Abel}

Let $L_{\infty }$ denote the unique $C^{\infty }$ complex--line bundle on $X$
with 
\begin{equation*}
c_{1}\left( L_{\infty }\right) =d.
\end{equation*}
If $P$ is an effective divisor of degree $d$ on $X$, then there is a
holomorphic structure 
\begin{equation*}
\overline{\partial }_{P}:A_{X}^{0}\left( L_{\infty }\right) \rightarrow
A_{X}^{0,1}\left( L_{\infty }\right)
\end{equation*}
on $L_{\infty }$ and a section $s_{P}$ of $L_{\infty }$ such that 
\begin{eqnarray*}
\overline{\partial }_{P}\left( s_{P}\right) &=&0 \\
div\left( s_{P}\right) &=&P.
\end{eqnarray*}
Given any metric $\mu $ on $L_{\infty }$, there is a unique metric-$\left(
1,0\right) $ connection $D_{\mu ,P}$ on $L_{\infty }$ such that 
\begin{equation*}
D_{\mu ,P}^{0,1}=\overline{\partial }_{P}.
\end{equation*}
Alternatively, if we restrict our consideration to $X-\left| P\right| $,
there is a unique connection $D_{P}$ on 
\begin{equation*}
\left. L_{\infty }\right| _{X-\left| P\right| }
\end{equation*}
such that $s_{P}$ is flat. Again for this connection 
\begin{equation*}
D_{P}^{0,1}=\overline{\partial }_{P}.
\end{equation*}

If $Q$ is another effective divisor of degree $d$, we consider both $s_{P}$
and $s_{Q}$ as $C^{\infty }$-sections of the $C^{\infty }$-line bundle $%
L_{\infty }$ and define as above the map 
\begin{equation*}
g=\frac{s_{Q}}{s_{P}}:X-\left( \left| P\right| \cup \left| Q\right| \right)
\rightarrow \Bbb{C}^{*}
\end{equation*}
where we assume that $g$ is meromorphic in a small analytic neighborhood of $%
\left( \left| P\right| \cup \left| Q\right| \right) $. We compute 
\begin{eqnarray*}
0 &=&D_{Q}\left( g\cdot s_{P}\right) =dg\cdot s_{P}+g\cdot D_{Q}s_{P} \\
D_{Q}s_{P} &=&-g^{-1}\cdot dg\cdot s_{P}.
\end{eqnarray*}
So, for 
\begin{equation*}
\alpha _{PQ}=D_{Q}-D_{P},
\end{equation*}
we have 
\begin{eqnarray*}
\alpha \cdot s_{P} &=&\left( D_{Q}-D_{P}\right) s_{P} \\
&=&D_{Q}s_{P} \\
&=&-g^{-1}\cdot dg\cdot s_{Q}
\end{eqnarray*}
so that 
\begin{eqnarray*}
\alpha _{PQ} &=&-g^{-1}\cdot dg \\
\alpha _{PQ}^{\left( 0,1\right) } &=&\overline{\partial }_{Q}-\overline{%
\partial }_{P} \\
&=&-g^{-1}\cdot \overline{\partial }g.
\end{eqnarray*}

Now via residue and the fact that 
\begin{equation*}
\underset{r\rightarrow 0}{\lim }r\mathrm{\log }r=0,
\end{equation*}
we have cohomologous currents 
\begin{eqnarray*}
2\pi i\int\nolimits_{0}^{\infty } &\sim &\int\nolimits_{X}d\theta \wedge  \\
&\sim &\int\nolimits_{X}\left( z^{-1}\cdot dz\right) \wedge .
\end{eqnarray*}
on $\overline{\Bbb{C}^{*}}=\Bbb{P}^{1}$. We therefore obtain the equality 
\begin{equation}
\int\nolimits_{P}^{Q}\eta =\frac{1}{2\pi i}\int\nolimits_{X}\alpha
_{PQ}\wedge \eta   \label{abel2}
\end{equation}
for $\eta \in H^{1,0}\left( X\right) $ by pulling back the above cohomology
between currents via $g$. Also, since 
\begin{eqnarray*}
D_{P}^{0,1} &=&D_{\mu ,P}^{0,1} \\
D_{Q}^{0,1} &=&D_{\mu ,Q}^{0,1}
\end{eqnarray*}
we have for 
\begin{equation*}
\alpha _{\mu ,PQ}:=D_{\mu ,Q}-D_{\mu ,P}
\end{equation*}
that 
\begin{equation*}
\alpha _{\mu ,PQ}^{0,1}=\alpha _{PQ}^{0,1}
\end{equation*}
so that the $\left( 0,1\right) $-summand of $-g^{-1}\cdot dg$ is bounded.

The equality $\left( \ref{abel2}\right) $ shows that if, for some 
\begin{equation*}
\varepsilon \in H_{1}\left( X;\Bbb{Z}\right) ,
\end{equation*}
one has 
\begin{equation*}
\int\nolimits_{P}^{Q}=\int\nolimits_{\varepsilon }:H^{1,0}\left( X\right)
\rightarrow \Bbb{C},
\end{equation*}
then the deRham class 
\begin{equation*}
A=A^{1,0}\left( with\ poles\right) +A^{0,1}\left( bounded\right)
\end{equation*}
such that 
\begin{equation*}
\left\{ A^{0,1}\right\} \in 2\pi i\cdot H^{1}\left( X;\Bbb{Z}\right)
+H^{1,0}\left( X\right) .
\end{equation*}
So the Poincar\'{e} dual of $\varepsilon $ is a deRham class 
\begin{equation*}
\left\{ \xi \right\} \in H^{1}\left( X;\Bbb{Z}\right)
\end{equation*}
such that 
\begin{equation*}
A^{0,1}-\xi ^{0,1}=\overline{\partial }\gamma
\end{equation*}
for some $C^{\infty }$-function 
\begin{equation*}
\gamma \in A_{X}^{0}.
\end{equation*}
Then the form 
\begin{equation*}
\psi :=g^{-1}\cdot dg+\xi +d\gamma
\end{equation*}
on $X-\left( \left| P\right| \cup \left| Q\right| \right) $ is $d$-closed
and of type $\left( 1,0\right) $ and therefore meromorphic on $X$ with poles
at the divisor $Q-P$.

As mentioned above, the ``second part'' of of the proof of Abel's theorem
consists in defining the meromorphic function 
\begin{equation*}
f=e^{\int \psi }
\end{equation*}
and noticing that 
\begin{equation*}
div\left( f\right) =Q-P.
\end{equation*}
Thus 
\begin{equation*}
f:X\rightarrow \Bbb{P}^{1}
\end{equation*}
gives the rational equivalence between $P$ and $Q$.

\subsection{Analogue for one-cycles on threefolds}

In the remainder of this paper we present the analogue of Abel's theorem in
the case in which $P$ and $Q$ are cohomologous (sums of) smooth
sub-canonical curves on a threefold $X$ and $E_{\infty }$ is a $C^{\infty }$
rank-$2$ vector bundle on $X$ whose first Chern class is trivial and whose
second Chern class is represented by $P$ (and so also by $Q$).

\section{Quaternionic connections}

\subsection{Serre's construction}

Let $X$ be a smooth projective threefold with 
\begin{equation*}
H^{1}\left( \mathcal{O}_{X}\right) =H^{2}\left( \mathcal{O}_{X}\right) =0
\end{equation*}
and that $P$ is a smooth (but possibly reducible) Riemann surface lying
inside $X$. We wish to consider situations in which 
\begin{equation}
\left. \omega _{X}\right| _{P}=\omega _{P},  \label{ass1}
\end{equation}
that is, the line bundle 
\begin{equation*}
\omega _{X}^{-1}\otimes \omega _{P}=\omega _{X}^{-1}\otimes
ext_{X}^{2}\left( \mathcal{O}_{P},\omega _{X}\right) =\omega
_{X}^{-1}\otimes ext_{X}^{2}\left( \mathcal{O}_{P},\mathcal{O}_{P}\right)
\end{equation*}
on $P$ is the trivial bundle$.$ Of course this is often not the case.
However if we assume, more generally, that there exists an \textit{effective}
divisor 
\begin{equation*}
\sum\nolimits_{j}x_{j}
\end{equation*}
on $P$ representing the line bundle 
\begin{equation*}
\omega _{X}^{-1}\otimes \omega _{P}
\end{equation*}
with $x_{j}$ distinct, and denote by 
\begin{equation}
\tilde{X}\rightarrow X  \label{blowdown}
\end{equation}
the blow-up of $X$ at the points $x_{j}$ containing the proper transform $%
\tilde{P}$ of $P$, then 
\begin{equation*}
\left. \omega _{\tilde{X}}\right| _{\tilde{P}}=\omega _{\tilde{P}}.
\end{equation*}

Assuming $\left( \ref{ass1}\right) $ from now on, we obtain 
\begin{eqnarray}
\mathcal{O}_{P} &=&ext_{X}^{2}\left( \mathcal{O}_{P},\mathcal{O}_{X}\right)
\label{critical} \\
&=&ext_{X}^{1}\left( \mathcal{I}_{P},\mathcal{O}_{X}\right) .  \notag
\end{eqnarray}
Since 
\begin{equation*}
H^{1}\left( hom_{X_{\delta }}\left( \mathcal{I}_{P},\mathcal{O}_{X}\right)
\right) =H^{1}\left( \mathcal{O}_{X}\right) =0,
\end{equation*}
we have 
\begin{equation*}
Ext^{1}\left( \mathcal{I}_{P},\mathcal{O}_{X}\right) =H^{0}\left(
ext_{X}^{1}\left( \mathcal{I}_{P},\mathcal{O}_{X}\right) \right)
\end{equation*}
and the nowhere-vanishing section of $ext_{X}^{1}\left( \mathcal{I}_{P},%
\mathcal{O}_{X}\right) $ gives, via the Serre construction, an exact
sequence of vector bundles 
\begin{equation}
0\rightarrow \mathcal{O}_{X}\rightarrow E\rightarrow \mathcal{I}
_{P}\rightarrow 0.  \label{defseq}
\end{equation}
Also $E$ has a distinguished global section 
\begin{equation*}
s_{P}
\end{equation*}
vanishing exactly at $P$. By the irreducibility of $P$, 
\begin{equation*}
H^{1}\left( \mathcal{I}_{P}\right) =0.
\end{equation*}
Thus 
\begin{equation*}
H^{1}\left( E\right) =0
\end{equation*}
and so $s$ deforms with each deformation of $\left( E,X\right) $. Also,
since 
\begin{equation*}
\det \left( E\right) =\mathcal{O}_{X}
\end{equation*}
we have 
\begin{equation}
E^{\vee }=E.  \label{dual}
\end{equation}

Tensor the surjection 
\begin{equation}
E^{\vee }\rightarrow \mathcal{I}_{P}  \label{surj}
\end{equation}
induced by $s$ with $\mathcal{O}_{P}$ to obtain a surjection of rank-$2$
bundles on $P$ which is therefore an isomorphism. Dualize to obtain an
isomorphism 
\begin{equation*}
N_{P\backslash X}\rightarrow \left. E\right| _{P}.
\end{equation*}
Thus we have the isomorphism 
\begin{equation}
H^{\cdot }\left( N_{P\backslash X}\right) =H^{\cdot +1}\left( \mathcal{I}
_{P}\otimes E\rightarrow E\right) .  \label{rel1}
\end{equation}

On the other hand, $\left( \ref{defseq}\right) $ and $\left( \ref{dual}
\right) $ give the exact sequence 
\begin{equation}
0\rightarrow E\rightarrow E^{\vee }\otimes E\rightarrow \mathcal{I}
_{P}\otimes E\rightarrow 0  \label{es2}
\end{equation}
from which comes 
\begin{equation}
H^{\cdot }\left( E^{\vee }\otimes E\right) \rightarrow H^{\cdot }\left( 
\mathcal{I}_{P}\otimes E\right) .  \label{rel2}
\end{equation}
Together $\left( \ref{rel1}\right) $ and $\left( \ref{rel2}\right) $ relate
the deformation functor of $E$ to that of $P$.

\subsection{Quaternionic line bundles}

Now 
\begin{equation*}
\det \left( E\right) =\mathcal{O}_{X}
\end{equation*}
allowing us to choose a non-vanishing holomorphic section 
\begin{equation*}
1\in H^{0}\left( \det \left( E\right) \right) .
\end{equation*}
We next chose a hermitian metric $\mu $ on $E$. We have an associated
structure of a \textit{quaternionic line bundle} on $E$ by defining 
\begin{equation*}
\left( \mathbf{j}\cdot s^{\prime }\right)
\end{equation*}
for any (locally defined) $C^{\infty }$-section $s^{\prime }$ of $E$ as the
unique element such that we have an equality of linear operators 
\begin{equation*}
\frac{s\wedge \left( \mathbf{j}\cdot s^{\prime }\right) }{1}=\mu \left(
s,s^{\prime }\right)
\end{equation*}
on sections $s$ of $E$. Thus 
\begin{eqnarray}
\left( \mathbf{j}\cdot s\right) &\bot &s  \label{quat} \\
s\wedge \left( \mathbf{j}\cdot s\right) &=&\left\| s\right\| ^{2}\cdot 1. 
\notag
\end{eqnarray}
Since 
\begin{equation*}
s\wedge \left( \mathbf{j}\cdot i\cdot s^{\prime }\right) =\mu \left(
s,i\cdot s^{\prime }\right) =-i\cdot \mu \left( s,s^{\prime }\right)
=s\wedge \left( -i\cdot \mathbf{j}\cdot s^{\prime }\right) ,
\end{equation*}
we have a well-defined (left) action on $E$ by the group $\Bbb{H}^{*}$ of
non-zero quaternions and the action of $\mathbf{j}\cdot $ is conjugate
linear. Since 
\begin{eqnarray*}
s\wedge \left( \mathbf{j}\cdot s\right) &=&-\left( \mathbf{j}\cdot \mathbf{j}%
\cdot s\right) \wedge \left( \mathbf{j}\cdot s\right) \\
&=&\left( \mathbf{j}\cdot s\right) \wedge \left( \mathbf{j}\cdot \mathbf{j}
\cdot s\right) \\
&=&\left\| \mathbf{j}\cdot s\right\| ^{2}\cdot 1
\end{eqnarray*}
$\mathbf{j}\cdot $ acts as an isometry. Also, for any two sections $s$ and $%
s^{\prime }$, we have that 
\begin{eqnarray*}
\mathbf{j}\cdot s\wedge \mathbf{j}\cdot s^{\prime } &=&\mu \left( \mathbf{j}
\cdot s,s^{\prime }\right) \\
&=&\overline{\mu \left( s^{\prime },\mathbf{j}\cdot s\right) } \\
&=&\overline{s^{\prime }\wedge \mathbf{j}\cdot \mathbf{j}\cdot s} \\
&=&\overline{s\wedge s^{\prime }}.
\end{eqnarray*}
For a fixed non-zero local section $s_{0}$ of $E$, we have the framing $%
\left( s_{0},\ \mathbf{j}\cdot s_{0}\right) $ as a complex vector bundle,
and, for any section $s$, we have the formula 
\begin{eqnarray}
s &=&\frac{\mu \left( s,s_{0}\right) }{\left\| s_{0}\right\| ^{2}}\cdot
s_{0}+\frac{\mu \left( s,\mathbf{j}\cdot s_{0}\right) }{\left\|
s_{0}\right\| ^{2}}\cdot \left( \mathbf{j}\cdot s_{0}\right)  \label{basis}
\\
&=&\frac{s\wedge \left( \mathbf{j}\cdot s_{0}\right) }{\left\| s_{0}\right\|
^{2}}\cdot s_{0}-\frac{s\wedge s_{0}}{\left\| s_{0}\right\| ^{2}}\cdot
\left( \mathbf{j}\cdot s_{0}\right)
\end{eqnarray}
with respect to the unitary basis 
\begin{equation*}
\left( \frac{s_{0}}{\left\| s_{0}\right\| },\ \mathbf{j}\cdot \frac{s_{0}}{%
\left\| s_{0}\right\| }\right) .
\end{equation*}

Now suppose that $D_{\mu ,P}$ denotes the metric $\left( 1,0\right) $
-connection on $E$ with respect to the metric $\mu $. Then we have that 
\begin{equation*}
d\left( \mu \left( s,s^{\prime }\right) \right) =\mu \left( D_{\mu .P}\left(
s\right) ,s^{\prime }\right) +\mu \left( s,D_{\mu .P}\left( s^{\prime
}\right) \right)
\end{equation*}
so that 
\begin{equation*}
d\left( s\wedge \mathbf{j}\cdot s^{\prime }\right) =D_{\mu ,P}\left(
s\right) \wedge \mathbf{j}\cdot s^{\prime }+s\wedge \mathbf{j}\cdot D_{\mu
,P}\left( s^{\prime }\right)
\end{equation*}
from which follows that 
\begin{equation*}
d\left( s\wedge s^{\prime }\right) =D_{\mu ,P}\left( s\right) \wedge
s^{\prime }-s\wedge \mathbf{j}\cdot D_{\mu ,P}\left( \mathbf{j}\cdot
s^{\prime }\right) .
\end{equation*}
Since 
\begin{equation*}
d\left( s\wedge s^{\prime }\right) =D_{\mu ,P}\left( s\right) \wedge
s^{\prime }+s\wedge D_{\mu ,P}\left( s^{\prime }\right)
\end{equation*}
it follows that $\mathbf{j}\cdot $ commutes with $D_{\mu ,P}$ so that $%
D_{\mu ,P}$ is a quaternionic connection. Since $-\mathbf{j}\cdot a\cdot 
\mathbf{j}=\overline{a}$, 
\begin{eqnarray*}
-\mathbf{j}\cdot D_{\mu ,P}^{1,0}\cdot \mathbf{j}\cdot &=&D_{\mu ,P}^{0,1} \\
-\mathbf{j}\cdot D_{\mu ,P}^{0,1}\cdot \mathbf{j}\cdot &=&D_{\mu ,P}^{1,0}.
\end{eqnarray*}
Thus 
\begin{equation*}
D_{\mu ,P}^{1,0}=-\mathbf{j}\cdot \overline{\partial }\cdot \mathbf{j}\cdot
\end{equation*}

Now let $s=s_{P}$, the holomorphic section of $E$. Our first goal is to
understand. If we write 
\begin{eqnarray*}
D_{\mu ,P}\left( s_{P}\right)  &=&D_{\mu ,P}^{1,0}\left( s_{P}\right) =-%
\mathbf{j}\cdot \overline{\partial }\cdot \mathbf{j}\cdot s_{P} \\
&=&\alpha _{P}\cdot s_{P}+\beta _{P}\cdot \mathbf{j}\cdot s_{P}
\end{eqnarray*}
for forms $\alpha _{P}$ and $\beta _{P}$ of type $\left( 1,0\right) $, then 
\begin{eqnarray*}
D_{\mu ,P}\left( \mathbf{j}\cdot s_{P}\right)  &=&\mathbf{j}\cdot D_{\mu
,P}s_{P}=\mathbf{j}\cdot \left( \alpha _{P}\cdot s_{P}+\beta _{P}\cdot 
\mathbf{j}\cdot s_{P}\right)  \\
&=&-\overline{\beta _{P}}\cdot s_{P}+\overline{\alpha _{P}}\cdot \mathbf{j}%
\cdot s_{P}.
\end{eqnarray*}
And since 
\begin{eqnarray*}
\left( D_{\mu ,P}\left( \mathbf{j}\cdot s_{P}\right) \right) ^{0,1}\wedge
s_{P} &=&\overline{\partial }\left( \mathbf{j}\cdot s_{P}\right) \wedge s_{P}
\\
&=&-\overline{\partial }\left( \left\| s_{P}\right\| ^{2}\right) 
\end{eqnarray*}
we can write 
\begin{equation}
D_{\mu ,P}\left( 
\begin{array}{l}
s_{P} \\ 
\mathbf{j}\cdot s_{P}
\end{array}
\right) =\left( 
\begin{array}{ll}
\partial \mathrm{\log }\left( \left\| s_{P}\right\| ^{2}\right)  & \beta _{P}
\\ 
-\overline{\beta _{P}} & \overline{\partial }\mathrm{\log }\left( \left\|
s_{P}\right\| ^{2}\right) 
\end{array}
\right) \left( 
\begin{array}{l}
s_{P} \\ 
\mathbf{j}\cdot s_{P}
\end{array}
\right) .  \label{matrix}
\end{equation}

With respect to the basis $\left( s_{P},\mathbf{j}\cdot s_{P}\right) $, the
curvature of the connection $D_{\mu ,P}$ then becomes 
\begin{equation*}
\left( 
\begin{array}{ll}
\overline{\partial }\partial \mathrm{\log }\left( \left\| s_{P}\right\|
^{2}\right)  & d\beta _{P} \\ 
-d\overline{\beta _{P}} & \partial \overline{\partial }\mathrm{\log }\left(
\left\| s_{P}\right\| ^{2}\right) 
\end{array}
\right) +\left( 
\begin{array}{ll}
-\beta _{P}\wedge \overline{\beta _{P}} & \left( \partial -\overline{%
\partial }\right) \mathrm{\log }\left( \left\| s_{P}\right\| ^{2}\right)
\wedge \beta _{P} \\ 
\overline{\beta _{P}}\wedge \left( \overline{\partial }-\partial \right) 
\mathrm{\log }\left( \left\| s_{P}\right\| ^{2}\right)  & \beta _{P}\wedge 
\overline{\beta _{P}}
\end{array}
\right) .
\end{equation*}
Since this matrix must be of type $\left( 1,1\right) $ we conclude 
\begin{equation}
\partial \beta _{P}=\beta _{P}\wedge \partial \mathrm{\log }\left( \left\|
s_{P}\right\| ^{2}\right)   \label{20part}
\end{equation}
and the curvature matrix becomes 
\begin{equation}
R_{\mu ,P}=\left( 
\begin{array}{ll}
\overline{\partial }\partial \mathrm{\log }\left( \left\| s_{P}\right\|
^{2}\right) -\beta _{P}\wedge \overline{\beta _{P}} & \overline{\partial }%
\beta _{P}+\beta _{P}\wedge \overline{\partial }\mathrm{\log }\left( \left\|
s_{P}\right\| ^{2}\right)  \\ 
-\partial \overline{\beta _{P}}-\overline{\beta _{P}}\wedge \partial \mathrm{%
\log }\left( \left\| s_{P}\right\| ^{2}\right)  & \partial \overline{%
\partial }\mathrm{\log }\left( \left\| s_{P}\right\| ^{2}\right) +\beta
_{P}\wedge \overline{\beta _{P}}
\end{array}
\right) .  \label{Pcurvature}
\end{equation}

\subsection{Chern-Simons functional}

Following \cite{DT} define the holomorphic Chern-Simons functional as
follows. For any two connections 
\begin{equation*}
D_{0},D_{1}
\end{equation*}
on a $C^{\infty }$-vector bundle $E$, form the connection 
\begin{equation*}
\tilde{D}:=D_{0}+t\left( D_{1}-D_{0}\right) +d_{t}
\end{equation*}
on 
\begin{equation*}
X\times \left[ 0,1\right]
\end{equation*}
and let 
\begin{equation*}
\tilde{R}=\tilde{D}^{2}
\end{equation*}
denote its curvature form. Let $A_{X}^{\cdot }$ denote the $\Bbb{C}$-deRham
complex in $X$ and let $F^{\cdot }$ denote the Hodge filtration on $%
A_{X}^{\cdot }$.

\begin{definition}
The \textit{holomorphic Chern-Simons current} $CS_{D_{0}}\left( D_{1}\right) 
$ is the current of type $\left( 1,2\right) +\left( 0,3\right) $ given by
the functional 
\begin{eqnarray*}
CS_{D_{0}}\left( D_{1}\right) :F^{2}A_{X}^{3}\rightarrow \Bbb{\ C} \\
\tau \mapsto \int\nolimits_{X\times \left[ 0,1\right] }\tau \wedge tr\left( 
\tilde{R}\wedge \tilde{R}\right) .
\end{eqnarray*}
(This is an extension of the standard definition of holomorphic Chern-Simons
functional which refers to the restriction of the above current to $%
F^{3}H^{3}\left( X\right) $.)
\end{definition}

Since $tr\left( \tilde{R}\wedge \tilde{R}\right) $ pulls back to an exact
form on the associated principal bundle of frames of $E$ (see \S 3 of \cite
{CS}), one has for three connections $D_{0},D_{1},D_{2}$ that 
\begin{equation}
CS_{D_{0}}\left( D_{1}\right) +CS_{D_{1}}\left( D_{2}\right)
=CS_{D_{0}}\left( D_{2}\right) .  \label{prinbun}
\end{equation}

Let 
\begin{equation*}
A:=D_{1}-D_{0}
\end{equation*}
We compute $\tilde{R}$ as follows 
\begin{equation*}
\left( D_{0}+tA+d_{t}\right) \circ \left( D_{0}+tA+d_{t}\right)
=R_{0}+t\cdot D_{0}A+t^{2}\cdot A\wedge A+dt\wedge A
\end{equation*}
so that 
\begin{equation*}
\tau \wedge \tilde{R}\wedge \tilde{R}=\tau \wedge dt\wedge 2A\wedge \left(
R_{0}+t\cdot D_{0}A+t^{2}\cdot A\wedge A\right)
\end{equation*}
and finally 
\begin{equation}
CS_{D_{0}}\left( D_{1}\right) \left( \tau \right) =\int\nolimits_{X}\tau
\wedge tr\left( A\wedge \left( R_{0}+D_{0}A+\frac{2}{3}A\wedge A\right)
\right) .  \label{CSdiff"}
\end{equation}

Suppose that $D_{0}^{0,1}$ and $D_{1}^{0,1}$ both give complex structures on 
$E$, that is 
\begin{equation*}
R_{0}^{0,2}=\left( D_{0}^{0,1}\right) ^{2}=0=\left( D_{1}^{0,1}\right) ^{2}.
\end{equation*}
Then 
\begin{eqnarray*}
0 &=&\left( D_{0}^{0,1}+A^{0,1}\right) ^{2} \\
&=&D_{0}^{0,1}A^{0,1}+A^{0,1}\wedge A^{0,1}.
\end{eqnarray*}
Thus in this case 
\begin{equation}
tr\left( A\wedge \left( R_{0}+D_{0}A+\frac{2}{3}A\wedge A\right) \right)
^{0,3}=-\frac{1}{3}tr\left( \left( A^{0,1}\right) ^{\wedge 3}\right) .
\label{03}
\end{equation}

On the other hand, suppose that $D_{0}$ and $D_{1}$ are both flat. Then $%
R_{0}=0$ and 
\begin{eqnarray*}
0 &=&\left( D_{0}+A\right) ^{2} \\
&=&D_{0}A+A\wedge A.
\end{eqnarray*}
so that $\left( \ref{CSdiff"}\right) $ becomes 
\begin{equation}
CS_{D_{0}}\left( D_{1}\right) \left( \tau \right) =-\frac{1}{3}%
\int\nolimits_{X}\tau \wedge tr\left( A\wedge A\wedge A\right) .
\label{CSflat}
\end{equation}

\subsection{Comparing connections via Chern-Simons theory}

Let $D_{P}$ denote the unique $\Bbb{H}$-connection on the restriction $%
E^{\prime }$ of $E$ to 
\begin{equation*}
X^{\prime }:=X-\left| P\right| 
\end{equation*}
such that 
\begin{equation*}
D_{P}\left( s_{P}\right) =0.
\end{equation*}
($D_{P}$ is of course flat.) We first apply the Chern-Simons theory to the
two connections 
\begin{eqnarray*}
D_{0} &=&D_{P} \\
D_{1} &=&D_{0,P}^{1,0}+D_{\mu ,P}^{0,1}=:D_{P}^{\prime }
\end{eqnarray*}
where, as above, $D_{\mu ,P}$ is the metric $\left( 1,0\right) $-connection
associated to $\mu $ and the complex structure on $E$. Thus $D_{P}^{\prime }$
is a $\left( 1,0\right) $-connection for the complex structure on $E$. Let 
\begin{equation*}
A_{P}:=D_{P}^{\prime }-D_{P}.
\end{equation*}
Then by $\left( \ref{matrix}\right) $%
\begin{equation*}
A_{P}\left( 
\begin{array}{l}
s_{P} \\ 
\mathbf{j}\cdot s_{P}
\end{array}
\right) =\left( 
\begin{array}{ll}
0 & 0 \\ 
-\overline{\beta _{P}} & \overline{\partial }\mathrm{\log }\left\|
s_{P}\right\| ^{2}
\end{array}
\right) \left( 
\begin{array}{l}
s_{P} \\ 
\mathbf{j}\cdot s_{P}
\end{array}
\right) .
\end{equation*}
Also, in terms of the basis $\left( s_{P},\mathbf{j}\cdot s_{P}\right) $ we
can compute the curvature 
\begin{eqnarray}
R_{P}^{\prime } &=&\left( D_{P}^{\prime }\right) ^{2} \\
&=&\left( D_{P}+A_{P}\right) ^{2}  \notag \\
&=&\left( 
\begin{array}{ll}
0 & 0 \\ 
-d\overline{\beta _{P}}-\overline{\partial }\mathrm{\log }\left\|
s_{P}\right\| ^{2}\wedge \overline{\beta _{P}} & \partial \overline{\partial 
}\mathrm{\log }\left\| s_{P}\right\| ^{2}
\end{array}
\right)   \label{curv'} \\
&=&\left( 
\begin{array}{ll}
0 & 0 \\ 
-\partial \overline{\beta _{P}} & \partial \overline{\partial }\mathrm{\log }%
\left\| s_{P}\right\| ^{2}
\end{array}
\right)   \notag
\end{eqnarray}
using $\left( \ref{20part}\right) $.

By $\left( \ref{CSdiff"}\right) $ and $\left( \ref{Dprime}\right) $, for any 
$\left( 3,0\right) +\left( 2,1\right) $-form $\tau $, the Chern-Simons
functional $CS_{D_{P}}\left( D_{P}^{\prime }\right) \left( \tau \right) $ is
given by the expression 
\begin{equation}
\int\nolimits_{X^{\prime }}tr\left( \overline{\partial }\mathrm{\log }%
\left\| s_{P}\right\| ^{2}\wedge \partial \overline{\partial }\mathrm{\log }%
\left\| s_{P}\right\| ^{2}\right) \wedge \tau   \label{CSmakeholo}
\end{equation}
By $\left( \ref{CSdiff1}\right) $ and Stokes theorem, if $\tau $ is $d$
-closed, this reduces to 
\begin{eqnarray}
CS_{D_{P}}\left( D_{P}^{\prime }\right) \left( \tau \right) 
&=&\int\nolimits_{X^{\prime }}\left( \overline{\partial }\mathrm{\log }%
\left\| s_{P}\right\| ^{2}\wedge \partial \overline{\partial }\mathrm{\log }%
\left\| s_{P}\right\| ^{2}\right) \wedge \tau   \label{CScompare} \\
&=&\int\nolimits_{X^{\prime }}\left( d\mathrm{\log }\left\| s_{P}\right\|
^{2}\wedge \partial \overline{\partial }\mathrm{\log }\left\| s_{P}\right\|
^{2}\right) \wedge \tau   \notag \\
&=&\int\nolimits_{\partial X^{\prime }}\left( \mathrm{\log }\left\|
s_{P}\right\| ^{2}\wedge \partial \overline{\partial }\mathrm{\log }\left\|
s_{P}\right\| ^{2}\right) \wedge \tau .  \notag
\end{eqnarray}
where $\int\nolimits_{\partial X^{\prime }}$ is computed as 
\begin{equation*}
\underset{\varepsilon \rightarrow 0}{\lim }\int\nolimits_{B_{\varepsilon
}}\left( \mathrm{\log }\left\| s_{P}\right\| ^{2}\wedge \partial \overline{%
\partial }\mathrm{\log }\left\| s_{P}\right\| ^{2}\right) \wedge \tau .
\end{equation*}
where $B_{\varepsilon }$ is a tubular neighborhood of $P$ in $X$, of radius $%
\varepsilon $ in some fixed metric. The integrand is bounded by 
\begin{equation*}
const.\cdot \frac{\left( \mathrm{\log }(\varepsilon )\right) ^{2}}{%
\varepsilon ^{2}}.
\end{equation*}
So the above integral is bounded by a constant multiple of $\varepsilon
.\left( \mathrm{\log }(\varepsilon )\right) ^{2}$, which tends to zero as $%
\varepsilon \rightarrow 0$. Thus 
\begin{equation*}
CS_{D_{P}}\left( D_{P}^{\prime }\right) \sim 0.
\end{equation*}

Finally since 
\begin{equation*}
D_{\mu ,P}-D_{P}^{\prime }=A^{1,0}
\end{equation*}
is of type $\left( 1,0\right) $ we again use $\left( \ref{CSdiff"}\right) $
to conclude by type that for any form $\tau $ of type $\left( 3,0\right)
+\left( 2,1\right) $, 
\begin{eqnarray}
CS_{D_{P}^{\prime }}\left( D_{\mu ,P}\right) \left( \tau \right)
&=&\int\nolimits_{X}\tau \wedge tr\left( A^{1,0}\wedge \left( R_{P}^{\prime
}+dA^{1,0}+\frac{2}{3}A^{1,0}\wedge A^{1,0}\right) \right)
\label{CSmake more holo} \\
&=&0  \notag
\end{eqnarray}
since $R_{P}^{\prime }$ is of type $\left( 1,1\right) $. So by the
additivity formula $\left( \ref{prinbun}\right) $ we conclude 
\begin{equation}
CS_{D_{P}}\left( D_{\mu ,P}\right) \sim 0,  \label{delbarfinite}
\end{equation}
that is, the $\left( 1,2\right) +\left( 0,3\right) $ current $%
CS_{D_{P}}\left( D_{\mu ,P}\right) $ is $d$-exact in the sense of currents.

\section{Analogue of the first part of Abel's theorem}

\subsection{Changing the curve $P$}

Let $X$ and $P$ be as in the previous section. Suppose now we have another
irreducible smooth curve $Q\subseteq X$ with 
\begin{equation*}
\omega _{X}^{-1}\otimes \omega _{Q}=\mathcal{O}_{Q}.
\end{equation*}
We suppose there is a fixed $C^{\infty }$-vector bundle, which we call $%
E_{\infty }$ such that, as $C^{\infty }$-vector bundles, we have
isomorphisms 
\begin{equation}
E_{P}\leftrightarrow E_{\infty }\leftrightarrow E_{Q}.  \label{Cinftyiso}
\end{equation}

We use a fixed hermitian structure $\mu $ on $E_{\infty }$ and a fixed
section $1$ of $\det \left( E_{\infty }\right) $ as above to make $E_{\infty
}$, and therefore also $E_{P}$ and $E_{Q}$, into quaternionic line bundles
such that the above correspondences are $C^{\infty }$-isomorphisms of
quaterionic line bundles. Again considering $s_{P}$ and $s_{Q}$ as sections
of the quaternionic line bundle $E_{\infty }$ we have 
\begin{equation*}
s_{Q}=g\cdot s_{P}
\end{equation*}
where, for $X^{\prime }=X-\left( \left| P\right| \cup \left| Q\right|
\right) $, 
\begin{equation*}
g=\left( 
\begin{array}{ll}
a & b \\ 
-\overline{b} & \overline{a}
\end{array}
\right) :X^{\prime }\rightarrow \Bbb{H}^{*}.
\end{equation*}
That is 
\begin{equation*}
s_{Q}=a\cdot s_{P}+b\cdot \left( \mathbf{j}\cdot s_{P}\right) .
\end{equation*}
Then 
\begin{eqnarray*}
\mathbf{j}\cdot s_{Q} &=&\mathbf{j}\cdot a\cdot s_{P}+\mathbf{j}\cdot b\cdot
\left( \mathbf{j}\cdot s_{P}\right) \\
&=&-\overline{b}\cdot s_{P}+\overline{a}\cdot \left( \mathbf{j}\cdot
s_{P}\right)
\end{eqnarray*}
so that the matrix $g$ gives the expression for the $\Bbb{C}$-basis $\left(
s_{Q},\mathbf{j}\cdot s_{Q}\right) $ in terms of the $\Bbb{C}$-basis $\left(
s_{P},\mathbf{j}\cdot s_{P}\right) $, that is, 
\begin{equation*}
\left( 
\begin{array}{l}
s_{Q} \\ 
\mathbf{j}\cdot s_{Q}
\end{array}
\right) =g\cdot \left( 
\begin{array}{l}
s_{P} \\ 
\mathbf{j}\cdot s_{P}
\end{array}
\right) .
\end{equation*}

We are interested in comparing two connections on $E^{\prime }$, namely the $%
\Bbb{H}$-connection $D_{P}$ with flat section $s_{P}$ and the $\Bbb{H}$
-connection $D_{Q}$ with flat section $s_{Q}$. Again we compute 
\begin{eqnarray*}
0 &=&D_{Q}\left( g\cdot s_{P}\right) =dg\cdot s_{P}+g\cdot D_{Q}s_{P} \\
D_{Q}s_{P} &=&-g^{-1}\cdot dg\cdot s_{P}.
\end{eqnarray*}
So 
\begin{equation*}
\left( D_{Q}-D_{P}\right) \left( s_{P}\right) .
\end{equation*}
Since $\Bbb{H}$-connections are determined by their values on a single
section, then in terms of the basis $\left( s_{P},\mathbf{j}\cdot
s_{P}\right) $, we have 
\begin{equation}
A_{PQ}:=D_{Q}-D_{P}=-g^{-1}\cdot dg.  \notag
\end{equation}

Since $D_{P}$ and $D_{Q}$ are both flat, 
\begin{eqnarray*}
0 &=&\left( D_{P}+A_{PQ}\right) ^{2} \\
&=&dA_{PQ}+A_{PQ}\wedge A_{PQ}.
\end{eqnarray*}
Write 
\begin{equation*}
h=\left[ 
\begin{array}{cc}
a & b \\ 
-\overline{b} & \overline{a}
\end{array}
\right] =r\cdot \left[ 
\begin{array}{cc}
u & v \\ 
-\overline{v} & \overline{u}
\end{array}
\right] 
\end{equation*}
with $u\overline{u}+v\overline{v}=1$. So for 
\begin{equation*}
\varkappa =\left[ 
\begin{array}{cc}
u & v \\ 
-\overline{v} & \overline{u}
\end{array}
\right] =u+v\cdot \mathbf{j}
\end{equation*}
we have 
\begin{equation*}
\varkappa ^{-1}=\overline{u}-v\cdot \mathbf{j}
\end{equation*}
since 
\begin{eqnarray*}
\left( u+v\cdot \mathbf{j}\right) \cdot \left( \overline{u}-v\cdot \mathbf{j}%
\right)  &=&u\overline{u}+v\cdot \mathbf{j}\cdot \overline{u}-uv\cdot 
\mathbf{j-}v\cdot \mathbf{j}\cdot v\cdot \mathbf{j} \\
&=&u\overline{u}+vu\cdot \mathbf{j}-uv\cdot \mathbf{j}+v\overline{v} \\
&=&1.
\end{eqnarray*}
The form 
\begin{eqnarray*}
h^{-1}\cdot dh &=&d\mathrm{\log }r+\overline{\varkappa }^{t}d\varkappa  \\
&=&r^{-1}dr+\left( \overline{u}-v\cdot \mathbf{j}\right) \left( du+dv\cdot 
\mathbf{j}\right)  \\
&=&r^{-1}dr+\left( \overline{u}du+vd\overline{v}\right) +\left( \overline{u}%
dv-vd\overline{u}\right) \cdot \mathbf{j}
\end{eqnarray*}
is a left-invariant $1$-form on 
\begin{equation*}
\Bbb{H}^{*}=\left( 0,+\infty \right) \times SU(2).
\end{equation*}
Since $u\overline{u}+v\overline{v}=1$, the form 
\begin{equation*}
\overline{u}du+vd\overline{v}=-\left( ud\overline{u}+vd\overline{v}\right) 
\end{equation*}
is purely imaginary.

Using the rule 
\begin{equation*}
\left( A+B\cdot \mathbf{j}\right) \wedge \left( C+D\cdot \mathbf{j}\right)
=\left( A\wedge C-B\wedge \overline{D}\right) +\left( A\wedge D+B\wedge 
\overline{C}\right) \cdot \mathbf{j}
\end{equation*}
we have that, if $\overline{A}=-A$, 
\begin{equation*}
\left( A+B\cdot \mathbf{j}\right) \wedge \left( A+B\cdot \mathbf{j}\right)
=-B\wedge \overline{B}-2\left( \overline{A}\wedge B\right) \cdot \mathbf{j}
\end{equation*}
and 
\begin{equation*}
\left( A+B\cdot \mathbf{j}\right) ^{\wedge ^{3}}=3\overline{A}\wedge B\wedge 
\overline{B}.
\end{equation*}
We compute 
\begin{eqnarray*}
\left( h^{-1}dh\right) ^{\wedge ^{3}} &=&r^{-1}dr\wedge \left( \left( 
\overline{u}du+vd\overline{v}\right) +\left( \overline{u}dv-vd\overline{u}%
\right) \cdot \mathbf{j}\right) ^{\wedge ^{2}} \\
&&+\left( \left( \overline{u}du+vd\overline{v}\right) +\left( \overline{u}%
dv-vd\overline{u}\right) \cdot \mathbf{j}\right) ^{\wedge ^{3}}
\end{eqnarray*}
and 
\begin{equation*}
\overline{\left( \overline{u}du+vd\overline{v}\right) }\wedge \left( 
\overline{u}dv-vd\overline{u}\right) =d\overline{u}dv.
\end{equation*}
So 
\begin{eqnarray*}
\left( \left( \overline{u}du+vd\overline{v}\right) +\left( \overline{u}dv-vd%
\overline{u}\right) \cdot \mathbf{j}\right) ^{\wedge ^{3}} &=&3d\overline{u}%
dv\wedge \left( ud\overline{v}-\overline{v}du\right) \\
&=&3\left( ud\overline{u}dvd\overline{v}-\overline{v}dud\overline{u}%
dv\right) .
\end{eqnarray*}

On the other hand 
\begin{equation*}
0=d\left( h\cdot h^{-1}\right) =\left( dh^{-1}\right) \cdot h+h^{-1}\cdot dh
\end{equation*}
so that 
\begin{equation}
d\left( h^{-1}\cdot dh\right) =-h^{-1}\cdot dh\wedge h^{-1}\cdot dh
\label{eq1}
\end{equation}
and 
\begin{equation}
\left( h^{-1}\cdot dh\right) ^{\wedge 3}=-d\mathrm{\log }r\wedge d\left(
\varkappa ^{-1}\cdot d\varkappa \right) +\left( \varkappa ^{-1}\cdot %
d\varkappa \right) ^{\wedge 3}.  \label{eq2}
\end{equation}
So 
\begin{equation}
tr\left( \left( \varkappa ^{-1}\cdot d\varkappa \right) ^{\wedge 3}\right)
=3\left( ud\overline{u}dvd\overline{v}-\overline{v}dud\overline{u}dv\right)
+3\overline{\left( ud\overline{u}dvd\overline{v}-\overline{v}dud\overline{u}%
dv\right) }  \label{eq3}
\end{equation}
is real and $d$-closed on $X^{\prime }$ as is 
\begin{equation*}
tr\left( \left( h^{-1}\cdot dh\right) ^{\wedge 3}\right) .
\end{equation*}
Since the left-invariant form $\left( \ref{eq3}\right) $ evaluated at $%
\left( u,v\right) =\left( 1,0\right) $ is 
\begin{equation*}
12d\left( \mathrm{Im}a\right) \wedge d\left( \mathrm{Re}b\right) \wedge
d\left( \mathrm{Im}b\right) 
\end{equation*}
we have 
\begin{eqnarray*}
\int\nolimits_{S^{3}}tr\left( \left( h^{-1}\cdot dh\right) ^{\wedge
3}\right)  &=&12\cdot vol\left( S^{3}\right)  \\
&=&24\pi ^{2}.
\end{eqnarray*}

We compare these connections on $\left. E_{\infty }\right| _{X^{\prime }}$
via the Chern-Simons functional. Now 
\begin{equation*}
A_{PQ}=-g^{*}\left( h^{-1}\cdot dh\right) 
\end{equation*}
for the invariant one-form $h^{-1}\cdot dh$ on $\Bbb{H}^{*}$. Then by $%
\left( \ref{eq1}\right) $ and $\left( \ref{eq2}\right) $ we have the formula 
\begin{eqnarray}
CS_{D_{P}}\left( D_{Q}\right) \left( \tau \right) 
&=&\int\nolimits_{X}g^{*}tr\left( h^{-1}\cdot dh\wedge d\left( h^{-1}\cdot %
dh\right) -\frac{2}{3}\left( h^{-1}\cdot dh\right) ^{\wedge 3}\right) \wedge
\tau   \label{preabel} \\
&=&\frac{1}{3}\int\nolimits_{X}g^{*}tr\left( \left( h^{-1}\cdot dh\right)
^{\wedge 3}\right) \wedge \tau .  \notag
\end{eqnarray}
Analogously to the classical Abel theorem on curves, via residue and the
fact that 
\begin{equation*}
\underset{r\rightarrow 0}{\lim }r\mathrm{\log }r=0,
\end{equation*}
we have cohomologous currents 
\begin{eqnarray*}
tr\left( \left( h^{-1}\cdot dh\right) ^{\wedge 3}\right)  &\sim &tr\left(
\varkappa ^{-1}\cdot d\varkappa \right) ^{\wedge 3} \\
&\sim &24\pi ^{2}\int\nolimits_{\left( 0,+\infty \right) }
\end{eqnarray*}
on $\overline{\Bbb{H}^{*}}=\Bbb{HP}^{1}$. Pulling back via $g$ we therefore
have by $\left( \ref{preabel}\right) $ that, whenever $\tau $ is a $d$%
-closed form of type $\left( 3,0\right) +\left( 2,1\right) $, 
\begin{equation}
CS_{D_{P}}\left( D_{Q}\right) \left( \tau \right) =\frac{24\pi ^{2}}{3}%
\int\nolimits_{P}^{Q}\tau .  \label{abel}
\end{equation}
Here on the right-hand side we integrate over the $3$-chain 
\begin{equation*}
g^{-1}\left( \left( 0,\infty \right) \right) ,
\end{equation*}
which indeed bounds $Q-P$. This is of course completely analogous, for rank-$
2$ vector bundles, to the main step in the above version of the proof of
classical Abel's theorem for line bundles on curves. Thus 
\begin{equation*}
CS_{D_{P}}\left( D_{Q}\right) =8\pi ^{2}\int\nolimits_{P}^{Q}
\end{equation*}
is $d$-exact in the sense of currents. Since 
\begin{equation*}
tr\left( \left( h^{-1}\cdot dh\right) ^{\wedge 3}\right) 
\end{equation*}
is $d$-closed on $\Bbb{H}^{*}$, $CS_{D_{P}}\left( D_{Q}\right) $ is the $%
\left( 1,2\right) +\left( 0,3\right) $ summand of a $d$-closed form on $%
X^{\prime }$, not of a $d$-closed current on $X$. Indeed, by what we have
just shown, 
\begin{eqnarray*}
\frac{1}{3}\int\nolimits_{X}g^{*}tr\left( \left( h^{-1}\cdot dh\right)
^{\wedge 3}\right) \wedge d\beta  &=&8\pi ^{2}\int\nolimits_{P}^{Q}d\beta  \\
&=&8\pi ^{2}\left( \int\nolimits_{Q}\beta -\int\nolimits_{P}\beta \right) .
\end{eqnarray*}

Finally by $\left( \ref{abel}\right) $, $\left( \ref{delbarfinite}\right) $
and the additivity property $\left( \ref{prinbun}\right) $, we have that 
\begin{equation*}
CS_{D_{\mu ,P}}\left( D_{\mu ,Q}\right) -CS_{D_{P}}\left( D_{Q}\right)
\end{equation*}
is a current coboundary and so, for any $d$-closed form $\tau $ of type $%
\left( 3,0\right) +\left( 2,1\right) $ that 
\begin{equation}
CS_{D_{\mu ,P}}\left( D_{\mu ,Q}\right) \left( \tau \right) =8\pi
^{2}\int\nolimits_{Q}^{P}\tau .  \label{abel"}
\end{equation}

\subsection{$P$ and $Q$ Abel-Jacobi equivalent}

Let 
\begin{equation*}
A_{\mu ,PQ}=D_{\mu ,Q}-D_{\mu ,P}
\end{equation*}
where we recall that 
\begin{eqnarray*}
A_{\mu ,PQ}^{0,1} &=&\overline{\partial }_{Q}-\overline{\partial }_{P} \\
A_{\mu ,PQ}^{1,0} &=&-\mathbf{j}\cdot \left( \overline{\partial }_{Q}-%
\overline{\partial }_{P}\right) \cdot \mathbf{j}.
\end{eqnarray*}
Let 
\begin{equation*}
\alpha _{\mu ,PQ}:=tr\left( A_{\mu ,PQ}\wedge \left( R_{\mu ,P}+dA_{\mu ,PQ}+%
\frac{2}{3}A_{\mu ,PQ}\wedge A_{\mu ,PQ}\right) \right)
\end{equation*}
be the $C^{\infty }$-deRham form giving $CS_{D_{\mu ,P}}\left( D_{\mu
,Q}\right) $ so that 
\begin{eqnarray*}
CS_{D_{\mu ,P}}\left( D_{\mu ,Q}\right) \left( \tau \right)
&=&\int\nolimits_{X}\tau \wedge \alpha _{\mu ,PQ} \\
&=&\int\nolimits_{X}\tau \wedge \alpha _{\mu ,PQ}^{\left( 1,2\right) +\left(
0,3\right) }.
\end{eqnarray*}

If $P$ and $Q$ are algebraically equivalent, we can arrange that $\Gamma $
with 
\begin{equation*}
\partial \Gamma =Q-P
\end{equation*}
be chosen to lie inside a (possibly reducible) algebraic surface on $X$.
Suppose now that 
\begin{equation*}
\tau =\overline{\partial }\beta ^{2,0}.
\end{equation*}
Then, since any $\left( 3,0\right) $-form restricts to zero on $\Gamma $, we
conclude by $\left( \ref{abel"}\right) $ that 
\begin{eqnarray*}
\int\nolimits_{X}\tau \wedge \alpha _{\mu ,PQ}^{\left( 1,2\right) +\left(
0,3\right) } &=&\int\nolimits_{X}\overline{\partial }\beta ^{2,0}\wedge
\alpha _{\mu ,PQ} \\
&=&8\pi ^{2}\int\nolimits_{\Gamma }\overline{\partial }\beta ^{2,0} \\
&=&8\pi ^{2}\int\nolimits_{\Gamma }d\beta ^{2,0} \\
&=&8\pi ^{2}\left( \int\nolimits_{Q}\beta ^{2,0}-\int\nolimits_{P}\beta
^{2,0}\right) \\
&=&0.
\end{eqnarray*}
Thus:

\begin{lemma}
If $P$ and $Q$ are algebraically equivalent, $\alpha _{\mu ,PQ}^{\left(
1,2\right) +\left( 0,3\right) }$ is $\overline{\partial }$-closed.
\end{lemma}

Referring to $\left( \ref{apq}\right) $ and writing the distribution-valued
differential 
\begin{eqnarray*}
\alpha _{PQ} &:&=tr\left( A_{PQ}\wedge \left( dA_{PQ}+\frac{2}{3}
A_{PQ}\wedge A_{PQ}\right) \right) \\
&=&-\frac{1}{3}tr\left( A_{PQ}\wedge A_{PQ}\wedge A_{PQ}\right) ,
\end{eqnarray*}
then 
\begin{eqnarray*}
CS_{D_{P}}\left( D_{Q}\right) \left( \tau \right) &=&\int\nolimits_{X}\tau
\wedge \alpha _{PQ} \\
&=&\int\nolimits_{X}\tau \wedge \alpha _{PQ}^{\left( 1,2\right) +\left(
0,3\right) }.
\end{eqnarray*}
Again if $P$ is algebraically equivalent to $Q$, by $\left( \ref{abel}%
\right) $, $\Gamma $ can be chosen so that 
\begin{equation*}
\int\nolimits_{X}\tau \wedge \alpha _{PQ}^{\left( 1,2\right) +\left(
0,3\right) }=0
\end{equation*}
whenever $\tau \in F^{2}A_{X}^{3}$ is $\overline{\partial }$-exact. Thus:

\begin{lemma}
If $P$ is algebraically equivalent to $Q$, then $\alpha _{PQ}^{\left(
1,2\right) +\left( 0,3\right) }$ is $\overline{\partial }$-closed.
\end{lemma}

The equality $\left( \ref{abel"}\right) $ shows that, if $P$ is
algebraically equivalent to $Q$ and 
\begin{equation*}
\int\nolimits_{Q}^{P}=\int\nolimits_{\varepsilon }:F^{2}H^{3}\left( X\right)
\rightarrow \Bbb{C}
\end{equation*}
for some 
\begin{equation*}
\varepsilon \in H_{3}\left( X;\Bbb{Z}\right) ,
\end{equation*}
then we have the containment 
\begin{equation*}
\left\{ \alpha _{\mu ,PQ}^{\left( 1,2\right) +\left( 0,3\right) }\right\}
\in H^{3}\left( X;\Bbb{Z}\right) +F^{2}H^{3}\left( X\right)
\end{equation*}
of the Dolbeault class of $\alpha _{\mu ,PQ}^{\left( 1,2\right) +\left(
0,3\right) }$. So the Poincar\'{e} dual of $\varepsilon $ is a deRham class 
\begin{equation*}
\left\{ \xi \right\} \in H^{3}\left( X;\Bbb{Z}\right)
\end{equation*}
with the property that 
\begin{equation*}
\alpha _{\mu ,PQ}^{\left( 1,2\right) +\left( 0,3\right) }-\xi ^{\left(
1,2\right) +\left( 0,3\right) }=\overline{\partial }\gamma _{\mu }
\end{equation*}
for some $C^{\infty }$-form 
\begin{equation*}
\gamma _{\mu }\in A_{X}^{\left( 1,1\right) +\left( 0,2\right) }.
\end{equation*}
Then the form 
\begin{equation*}
\psi _{\mu }:=\alpha _{\mu ,PQ}-\xi -\overline{\partial }\gamma _{\mu }
\end{equation*}
on $X$ is of type $\left( 3,0\right) +\left( 2,1\right) $.

The equality $\left( \ref{abel}\right) $ shows that, if $P$ is algebraically
equivalent to $Q$ and 
\begin{equation*}
\int\nolimits_{Q}^{P}=\int\nolimits_{\varepsilon }:F^{2}H^{3}\left( X\right)
\rightarrow \Bbb{C}
\end{equation*}
for some 
\begin{equation*}
\varepsilon \in H_{3}\left( X;\Bbb{Z}\right) ,
\end{equation*}
then we have as above, taking the deRham dual 
\begin{equation*}
\left\{ \xi \right\} \in H^{3}\left( X;\Bbb{Z}\right)
\end{equation*}
of $\varepsilon $, we have 
\begin{equation*}
\int\nolimits_{X}\left( \alpha _{PQ}-\xi \right) \wedge \tau =0
\end{equation*}
for any $d$-closed form $\tau $ of type $\left( 3,0\right) +\left(
2,1\right) $ and any $\overline{\partial }$-exact form of type $\left(
3,0\right) +\left( 2,1\right) $, hence for any $\overline{\partial }$-closed
form of type $\left( 3,0\right) +\left( 2,1\right) $. This means that the
current $\left( \alpha _{PQ}-\xi \right) ^{\left( 1,2\right) +\left(
0,3\right) }$ is $\overline{\partial }$-exact, and so, in particular, 
\begin{equation*}
\alpha _{PQ}^{0,3}-\xi ^{0,3}=\overline{\partial }\gamma ^{0,2}
\end{equation*}
for some distribution--valued form $\gamma ^{0,2}$ of type $\left(
0,2\right) $. Then the current 
\begin{equation*}
\alpha _{PQ}-\xi -d\gamma ^{0,2}
\end{equation*}
on $X$ is of type $\left( 3,0\right) +\left( 2,1\right) +\left( 1,2\right) $
and is $d$-closed on $X^{\prime }$. Since $\alpha _{PQ}-\xi -d\gamma $
integrates to zero against any $\overline{\partial }$-closed form of type $%
\left( 3,0\right) +\left( 2,1\right) $ on $X$, the current $\left( \alpha
_{PQ}-\xi -d\gamma \right) ^{1,2}$ is again $\overline{\partial }$-exact, we
can find a current $\gamma ^{1,1}$ such that the distribution valued-form 
\begin{equation*}
\psi :\alpha _{PQ}-\xi -d\gamma ^{0,2}-d\gamma ^{1,1}
\end{equation*}
on $X$ is a current of type $\left( 3,0\right) +\left( 2,1\right) $ and is $%
d $-closed on $X^{\prime }$.

Of course, in the case of classical Abel's theorem, $P$ and $Q$ are always
algebraically equivalent and $\psi $ gives directly the rational equivalence
of $P$ and $Q$. $\psi _{\mu }$ does not enter the picture because $\alpha
^{0,1}$ ($=\alpha _{\mu }^{0,1}$ for any metric $\mu $) is already smooth.
In the threefold case, $\left( \ref{CSmakeholo}\right) $ and $\left( \ref
{CSmake more holo}\right) $ imply that 
\begin{equation*}
\alpha _{PQ}^{0,3}=\alpha _{\mu ,PQ}^{0,3}
\end{equation*}
for any metric $\mu $ but in general the currents $\alpha _{PQ}^{1,2}$ and $%
\alpha _{\mu ,PQ}^{1,2}$ are not equal. However it should be true that, as
we vary the metric $\mu $ nicely so that it becomes flat on $X^{\prime }$,
we achieve that $\alpha _{\mu ,PQ}^{1,2}$ converges to $\alpha _{PQ}^{1,2}$.

\end{document}